\documentclass[onecolumn,11pt]{article}
\usepackage[top=1in, bottom=1in, left=1in, right=1in]{geometry}
\usepackage{amsmath,amsfonts,amscd,amssymb}
\usepackage{graphicx}
\usepackage{epstopdf}
\usepackage{overpic}
\usepackage{cancel}
\usepackage{rotating}
\usepackage{url}
\usepackage{caption}
\usepackage{color}
\usepackage{rotating}
\usepackage{multirow}
\usepackage{wrapfig}
\usepackage{mathtools}
\usepackage{subeqnarray}
\usepackage{setspace}
\usepackage{palatino} 
\usepackage[numbers,sort&compress]{natbib}

\usepackage[bottom,flushmargin,hang,multiple]{footmisc}
\usepackage{lipsum}
\newcommand\blfootnote[1]{%
  \begingroup
  \renewcommand\thefootnote{}\footnote{#1}%
  \addtocounter{footnote}{-1}%
  \endgroup
}

\usepackage{hyperref}
\hypersetup{
    colorlinks=true,
    linkcolor=blue,
    filecolor=magenta,      
    urlcolor=blue,
    citecolor=black,
    }

\definecolor{header1}{cmyk}{0,0,0,1}

\DeclareMathOperator*{\argmin}{arg\rm{}min}

\usepackage{listings}
\usepackage{xcolor}

\definecolor{codegreen}{rgb}{0.0,0.3,0.0}
\definecolor{codegray}{rgb}{0.5,0.5,0.5}
\definecolor{codepurple}{rgb}{0.1,0,0.82}
\definecolor{backcolour}{rgb}{1,1,1}

\lstdefinestyle{mystyle}{
    backgroundcolor=\color{backcolour},   
    commentstyle=\color{codegreen},
    keywordstyle=\bf\color{black},
    numberstyle=\tiny\color{codegray},
    stringstyle=\color{codepurple},
    basicstyle=\ttfamily\small,
    breakatwhitespace=true,         
    breaklines=true,                 
    captionpos=t,                    
    keepspaces=true,                 
    numbers=none,
    frame=lines,
    numbersep=5pt,                  
    showspaces=false,                
    showstringspaces=false,
    showtabs=false,                  
    tabsize=3
}

\usepackage{mathabx}
\usepackage{subcaption}

\usepackage[utf8]{inputenc}

\title{\vspace{-.55in}{\fontsize{17.5}{17.5}\selectfont\textbf{Ensemble-SINDy: Robust sparse model discovery in the low-data, high-noise limit, with active learning and control}}}

\author{\normalsize{Urban Fasel$^{1*}$, J. Nathan Kutz$^2$, Bingni W. Brunton$^3$, Steven L. Brunton$^1$}\\
\footnotesize{$^1$ Department of Mechanical Engineering, University of Washington, USA}\\
\footnotesize{$^2$ Department of Applied Mathematics, University of Washington, USA}\\
\footnotesize{$^3$ Department of Biology, University of Washington, USA \vspace{-.1in}}
}

\date{}
\begin{document}
\maketitle

\lstdefinestyle{interfaces}{
  float,
  floatplacement=tbp
}

\blfootnote{$^*$ Corresponding author (ufasel@uw.edu).}
\vspace{-.2in}
\begin{abstract}
Sparse model identification enables the discovery of nonlinear dynamical systems purely from data; however, this approach is sensitive to noise, especially in the low-data limit.   
In this work, we leverage the statistical approach of bootstrap aggregating (bagging) to robustify the sparse identification of nonlinear dynamics (SINDy) algorithm. 
First, an ensemble of SINDy models is identified from subsets of limited and noisy data.  
The aggregate model statistics are then used to produce inclusion probabilities of the candidate functions, which enables uncertainty quantification and probabilistic forecasts.  
We apply this ensemble-SINDy (E-SINDy) algorithm to several synthetic and real-world data sets and demonstrate substantial improvements to the accuracy and robustness of model discovery from extremely noisy and limited data. 
For example, E-SINDy uncovers partial differential equations models from data with more than twice as much measurement noise as has been previously reported. 
Similarly, E-SINDy learns the Lotka Volterra dynamics from remarkably limited data of yearly lynx and hare pelts collected from 1900-1920.  
E-SINDy is computationally efficient, with similar scaling as standard SINDy.  
Finally, we show that ensemble statistics from E-SINDy can be exploited for active learning and improved model predictive control.

\end{abstract}

\section{Introduction}

Data-driven model discovery enables the characterization of complex systems where first principles derivations remain elusive, such as in neuroscience, power grids, epidemiology, finance, and ecology.
A wide range of data-driven model discovery methods exist, including equation-free modeling~\cite{kevrekidis2003equation}, normal form identification~\cite{majda2009normal,Yair2017pnas,kalia2021learning}, nonlinear Laplacian spectral analysis~\cite{giannakis2012nonlinear}, Koopman analysis~\cite{Mezic2005nd,Mezic2013arfm} and dynamic mode decomposition (DMD)~\cite{Schmid2010jfm,Rowley2009jfm,Kutz2016book}, symbolic regression~\cite{bongard2007automated, schmidt2009distilling,schmidt2011automated, daniels2015automated, daniels2015efficient}, sparse regression~\cite{Brunton2016pnas,rudy2017data}, Gaussian processes~\cite{raissi2017machine}, and deep learning~\cite{raissi2019physics,chen2018neural,Champion2019pnas,li2020fourier,rackauckas2020universal,lu2021learning}. 
Limited data and noisy measurements are fundamental challenges for all of these model discovery methods, often limiting the effectiveness of such techniques across diverse application areas.
The {\em sparse identification of nonlinear dynamics} (SINDy)~\cite{Brunton2016pnas} algorithm is promising, because it enables the discovery of interpretable and generalizable models that balance accuracy and efficiency.  
Moreover, SINDy is based on simple sparse linear regression that is highly extensible and requires significantly less data in comparison to, for instance, neural networks. 
In this work, we unify and extend innovations of the SINDy algorithm by leveraging classical statistical bagging methods~\cite{breiman1996bagging} to produce a computationally efficient and robust probabilistic model discovery method that overcomes the two canonical failure points of model discovery: noise and limited data.


The SINDy algorithm~\cite{Brunton2016pnas} provides a data-driven model discovery framework, relying on sparsity-promoting optimization to identify parsimonious models that avoid overfitting. 
These models may be ordinary differential equations (ODEs)~\cite{Brunton2016pnas} or partial differential equations (PDEs) \cite{rudy2017data,Schaeffer2017prsa}.
SINDy has been applied to a number of challenging model discovery problems, including for reduced-order models of fluid dynamics~\cite{Loiseau2017jfm,Loiseau2018jfm,loiseau2020data,guan2020sparse,deng2021galerkin,callaham2021empirical} and plasma dynamics~\cite{Dam2017pf,alves2020data,kaptanoglu2020physics}, turbulence closures~\cite{beetham2020formulating,beetham2021sparse,schmelzer2020discovery}, mesoscale ocean closures~\cite{zanna2020data}, nonlinear optics~\cite{Sorokina2016oe}, computational chemistry~\cite{boninsegna2018sparse}, and numerical integration schemes~\cite{Thaler2019jcp}.
SINDy has been widely adopted, in part, because it is highly extensible.  
Extensions of the SINDy algorithm include accounting for control inputs~\cite{Kaiser2018prsa} and rational functions~\cite{mangan2016inferring, kaheman2020sindy}, enforcing known conservation laws and symmetries~\cite{Loiseau2017jfm}, promoting stability~\cite{kaptanoglu2021promoting}, 
improved noise robustness through the integral formulation~\cite{schaeffer2017sparse,reinbold2020using,alves2020data, reinbold2021robust, messenger2021weak1, messenger2021weak2}, generalizations for stochastic dynamics~\cite{boninsegna2018sparse,callaham2021nonlinear} and
tensor formulations~\cite{Gelss2019mindy}, 
and probabilistic model discovery via sparse Bayesian inference~\cite{galioto2020bayesian, niven2020bayesian, yang2020bayesian, hirsh2021sparsifying,delahunt2021toolkit}.
Many of these innovations have been incorporated into the open source software package PySINDy~\cite{deSilva2020JOSS,kaptanoglu2021pysindy}.
Today, the biggest challenge with SINDy, and more broadly in model discovery, is learning models from limited and noisy data, especially for spatio-temporal systems governed by PDEs.


Model discovery algorithms are sensitive to noise because they rely on the accurate computation of derivatives, which is especially challenging for PDEs where noise can be strongly amplified when computing higher-order spatial derivatives.
There have been two key innovations to improve the noise robustness of SINDy:  control volume formulations and ensemble methods.  
The integral formulation of SINDy~\cite{schaeffer2017sparse} has proven powerful, enabling the identification of PDEs in a weak form that averages over control volumes, which significantly improves its noise tolerance. 
This approach has been used to discover a hierarchy of PDE models for fluids and plasmas~\cite{gurevich2019robust,reinbold2020using,alves2020data,gurevich2021learning,reinbold2021robust,messenger2021weak1,messenger2021weak2}. 
Several works have begun to explore ensemble methods to robustify data-driven modeling, including the use of bagging for DMD~\cite{sashidhar2021bagging}, ensemble-Lasso~\cite{sachdeva2019pyuoi}, subsample aggregating for improved discovery~\cite{nardini2020learning,delahunt2021toolkit}, statistical learning of PDEs to select model coefficients with high importance measures~\cite{maddu2019stability}, and improved discovery using ensembles based on sub-sampling of the data~\cite{reinbold2020using,gurevich2021learning,reinbold2021robust,delahunt2021toolkit}. 
When dealing with noise compromised data, it is also critical to provide uncertainty estimates of the discovered models. 
In this direction, recent innovations of SINDy use sparse Bayesian inference for probabilistic model discovery~\cite{galioto2020bayesian, niven2020bayesian, yang2020bayesian, hirsh2021sparsifying}.  Such methods employ Markov Chain Monte Carlo, which is extremely computationally intensive. 
These extensions have all improved the robustness of SINDy for high-noise data, although they have been developed largely in isolation and they have not been fully characterized, exploited, and/or integrated.  

\begin{figure*}[t]
    \includegraphics[width=\textwidth]{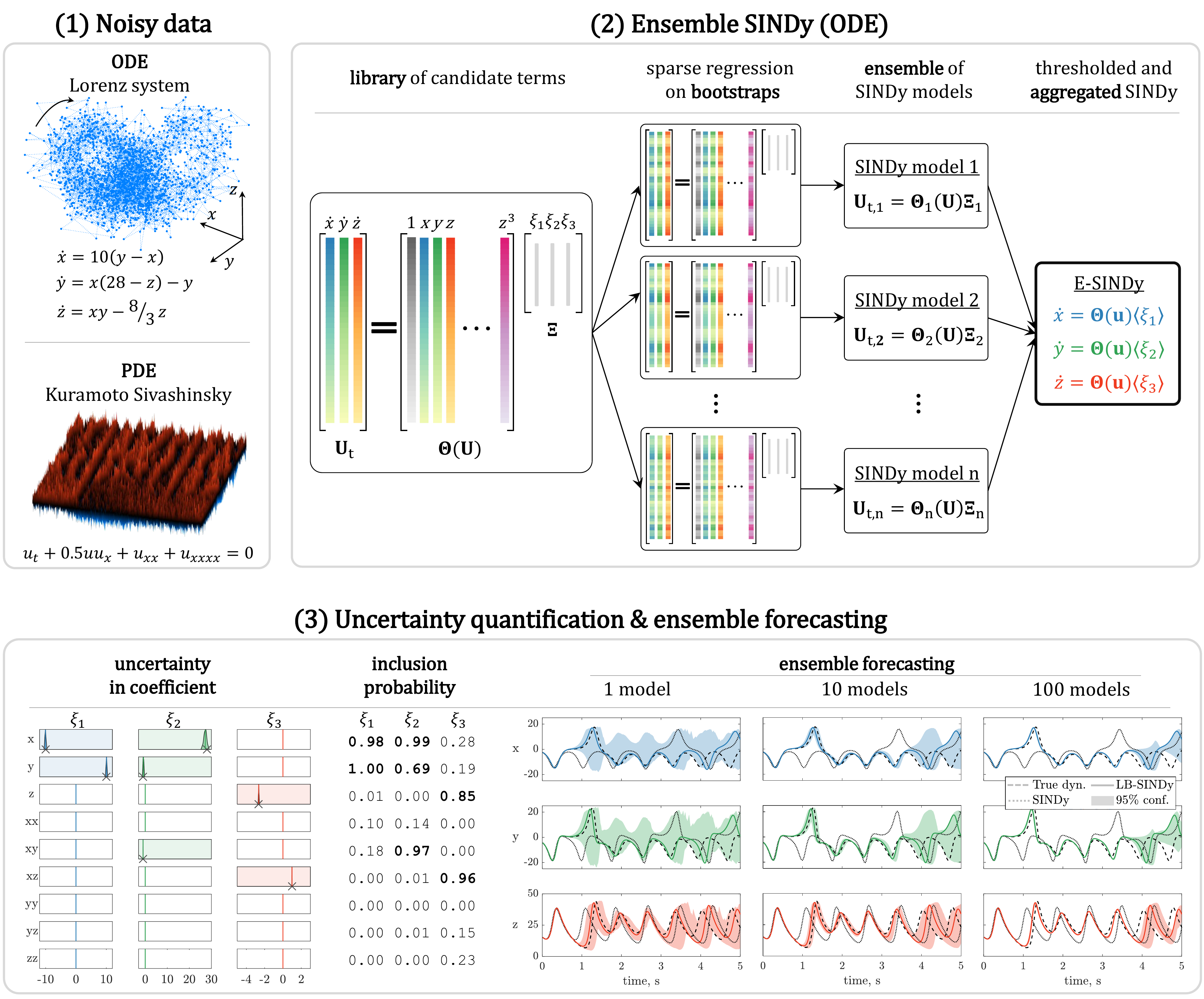}
    \vspace{-.275in}
    \caption{Schematic of the E-SINDy framework. E-SINDy exploits the statistical method of bootstrap aggregating (bagging) to identify ordinary and partial differential equations that govern the dynamics of observed noisy data. First, sparse regression is performed on bootstraps of the measured data, or on the library terms in case of library bagging, to identify an ensemble of SINDy models. The mean or median of the coefficients are then computed, coefficients with low inclusion probabilities are thresholded, and the E-SINDy model is aggregated and used for forecasting. }
    \label{fig1}
    \vspace{-.1in}
\end{figure*}

In this work, we unify and extend recent innovations in ensembling and the weak formulation of SINDy to develop and characterize a more robust E-SINDy algorithm.  
Further, we show how this method can be useful for active learning and control.
In particularly, we apply b(r)agging\footnote{b(r)agging refers to robust bagging based on median, rather than mean, information.} to SINDy to identify models of nonlinear ordinary differential equations of the form
\begin{equation}\label{eq1}
    \frac{\mathrm{d}}{\mathrm{dt}}\mathbf{u} = \mathbf{f}(\mathbf{u}), \hspace{10pt} \mathbf{u}(0) = \mathbf{u}_0,
\end{equation}
with state $\mathbf{u}\in\mathbb{R}^n$ and dynamics $\mathbf{f}(\mathbf{u})$
, and for nonlinear partial differential equations of the form
\begin{equation}\label{eq1b}
    \mathbf{u}_t = \mathbf{N}(\mathbf{u},\mathbf{u}_x,\mathbf{u}_{xx},\cdots,x,\mu),
\end{equation}
with $\mathbf{N}(\cdot)$ a system of nonlinear functions of the state $\mathbf{u}(x,t)$, its derivatives, and parameters $\mu$; partial derivatives are denoted with subscripts, such that $\mathbf{u}_t:= \partial \mathbf{u}/\partial t$.
We show that b(r)agging improves the accuracy and robustness of SINDy. 
The method also promotes interpretability through the inclusion  probabilities of candidate functions, enabling uncertainty quantification.  Importantly, the ensemble statistics are useful for producing probabilistic forecasts and can be used for active learning and nonlinear control.
We also demonstrate {library E-SINDy}, which subsamples terms in the SINDy library. 
E-SINDy is computationally efficient compared to probabilistic model identification methods based on Markov Chain Monte Carlo sampling~\cite{hirsh2021sparsifying}, which take several hours of CPU time to identify a model. 
In contrast, our method identifies models and summary statistics in seconds by leveraging the sparse regression of SINDy with statistical bagging techniques. 
Indeed, E-SINDy has similar computational scaling to standard SINDy.

We investigate different ensemble methods, apply them to several synthetic and real world data sets, and demonstrate that E-SINDy outperforms existing sparse regression methods, especially in the low-data and high-noise limit. 
A schematic of E-SINDy is shown in figure~\ref{fig1}.
We first describe SINDy for ordinary and partial differential equations in section~\ref{background}, before introducing the E-SINDy extension in section~\ref{method}, and discussing applications to challenging model discovery, active learning, and control problems in section~\ref{results}. 

\section{Background} \label{background}

Here we describe SINDy~\cite{Brunton2016pnas}, a data-driven model discovery method to identify sparse nonlinear models from measurement data. First, we introduce SINDy to identify ordinary differential equations, followed by its generalization to identify partial differential equations~\cite{rudy2017data,Schaeffer2017prsa}.

\subsection{Sparse identification of nonlinear dynamics}

The SINDy algorithm~\cite{Brunton2016pnas} identifies nonlinear dynamical systems from data, based on the assumption that many systems have relatively few active terms in the dynamics $\mathbf{f}$ in Eq.~\eqref{eq1}. SINDy uses sparse regression to identify these few active terms out of a library of candidate linear and nonlinear model terms.
Therefore, sparsity-promoting techniques may be used to find parsimonious models that automatically balance model complexity with accuracy~\cite{Brunton2016pnas}. 
We first measure $m$ snapshots of the state $\mathbf{u}$ in time and arrange these into a data matrix:
\begin{equation}
    \mathbf{U} = \left[\mathbf{u}_1 \hspace{3pt} \mathbf{u}_2 \hspace{3pt} \cdots \hspace{3pt} \mathbf{u}_m\right]^T.
\end{equation}
Next, we compute the library of $D$ candidate nonlinear functions $\boldsymbol{\Theta}(\mathbf{U})\in\mathbb{R}^{m\times D}$:
\begin{equation}
    \boldsymbol{\Theta}(\mathbf{U}) = [\mathbf{1} \hspace{6pt} \mathbf{U} \hspace{6pt} \mathbf{U}^2 \hspace{3pt} \cdots \hspace{3pt} \mathbf{U}^d \hspace{3pt} \cdots \hspace{3pt} \sin(\mathbf{U}) \hspace{3pt} \cdots ].
\end{equation} 
This library may include any functions that might describe the data, and this choice is crucial. The recommended strategy is to start with a basic choice, such as low-order polynomials, and then increase the complexity and order of the library until sparse and accurate models are obtained. 

We must also compute the time derivatives of the state $\mathbf{U}_t = \left[\dot{\mathbf{u}}_{1} \hspace{3pt} \dot{\mathbf{u}}_{2} \hspace{3pt} \cdots \hspace{3pt} \dot{\mathbf{u}}_{m}\right]^T$, typically by numerical differentiation.  
The system in Eq.~\eqref{eq1} may then be written in terms of these data matrices:
\begin{equation}
    \mathbf{U}_t = \boldsymbol{\Theta}(\mathbf{U})\boldsymbol{\Xi}.
\end{equation}
Each entry in $\boldsymbol{\Xi}\in\mathbb{R}^{D\times n}$ is a coefficient corresponding to a term in the dynamical system. Many dynamical systems have relatively few active terms in the governing equations. 
Thus, we may employ sparse regression to identify a sparse matrix of coefficients $\boldsymbol{\Xi}$ signifying the fewest active terms from the library that result in a good model fit: 
\begin{equation}
    \boldsymbol{\Xi} = \argmin_{\Hat{\boldsymbol{\Xi}}} \frac{1}{2} \|\mathbf{U}_t-\boldsymbol{\Theta}(\mathbf{U})\Hat{\boldsymbol{\Xi}}\|_2^2 + R(\hat{\boldsymbol{\Xi}}).
\end{equation}
The regularizer $R(\boldsymbol{\Xi})$ is chosen to promote sparsity in $\boldsymbol{\Xi}$.  
For example, sequentially thresholded least squares (STLS)~\cite{Brunton2016pnas} uses $R(\boldsymbol{\Xi})=\lambda \|\boldsymbol{\Xi}\|_0$ with a single hyperparameter $\lambda$, whereas sequentially thresholded ridge regression (STRidge)~\cite{rudy2017data} uses $R(\boldsymbol{\Xi})=\lambda_1 \|\boldsymbol{\Xi}\|_0 + \lambda_2 \|\boldsymbol{\Xi}\|_2$ with two hyperparameters $\lambda_1$ and $\lambda_2$. 
STLS was first introduced to discover ODEs and STRidge was introduced to discover PDEs where data can be highly correlated and STLS tends to perform poorly.
There are several other recently proposed regularizers and optimization schemes~\cite{champion2020unified,kaptanoglu2021promoting,carderera2021cindy}. 
We illustrate STRidge in pseudo code algorithm~\ref{algoSparsify}, noting that STRidge reduces to STLS for $\lambda_2=0$. 

\begin{lstlisting}[style=interfaces,language=Matlab,caption={Sparsify dynamics using STRidge (or STLS for $\lambda_2=0$).},label={algoSparsify},mathescape=true]
$\textbf{Inputs}:$ library $\boldsymbol{\Theta}$ and time derivatives $\mathbf{U}_t$
        sparsification and ridge parameters $\lambda_1$ and $\lambda_2$
$\textbf{Output}:$ SINDy model coefficients $\boldsymbol{\Xi}$

function sparsifyDynamics($\boldsymbol{\Theta}$, $\mathbf{U}_t$, $\lambda_1$, $\lambda_2$)
    for i = 1:n % n is state dimension
        % Initial guess: ridge regression
        $\boldsymbol{\Xi}_i$ = $(\boldsymbol{\Theta}^T_i \boldsymbol{\Theta}_i +\lambda_2 I)$ $\setminus$ $\boldsymbol{\Theta}^T_i\mathbf{U}_{t,i}$
        % Sequentially threshold and recompute ridge regression
        for k=1:10
            smallinds = (abs($\boldsymbol{\Xi}_i$)<$\lambda_1$); % Find small coefficients
            $\boldsymbol{\Xi}_i$(smallinds)=0; % and threshold
            bi = $\sim$smallinds(:,i) % big coefficients
            % Regress dynamics onto remaining terms to find sparse $\boldsymbol{\Xi}_i$
            $\boldsymbol{\Xi}_i$(bi) = $(\boldsymbol{\Theta}^T_i$(bi)$ \boldsymbol{\Theta}_i$(bi)$ +\lambda_2 I)$ $\setminus$ $\boldsymbol{\Theta}^T_i$(bi)$\mathbf{U}_{t,i}$
        end
        % Assemble identified model coefficients
        $\boldsymbol{\Xi}$(:,i) = $\boldsymbol{\Xi}_i$
    end
end
\end{lstlisting}

\subsection{Discovering partial differential equations}

 SINDy was recently generalized to identify partial differential equations~\cite{rudy2017data,Schaeffer2017prsa} in the \emph{partial differential equation functional identification of nonlinear dynamics} (PDE-FIND) algorithm. 
 PDE-FIND is similar to SINDy, but with the library including partial derivatives. 
 Spatial time-series data is arranged into a column vector $\mathbf{U}\in\mathbb{R}^{mn}$, with data collected over $m$ time points and $n$ spatial locations. Thus, for PDE-FIND, the library of candidate terms is $\boldsymbol{\Theta}(\mathbf{U})\in\mathbb{R}^{mn\times D}$. 
The PDE-FIND implementation of Rudy et al.~\cite{rudy2017data} takes derivatives using finite difference for clean data or polynomial interpolation for noisy data. The library of candidate terms can then be evaluated:
\begin{equation}
    \boldsymbol{\Theta}(\mathbf{U}) = [\mathbf{1} \hspace{6pt} \mathbf{U} \hspace{6pt} \mathbf{U}^2 \hspace{3pt} \cdots \hspace{3pt} \mathbf{U}_x \hspace{3pt} \cdots \hspace{3pt} \mathbf{UU}_x \hspace{3pt} \cdots ].
\end{equation} 
The time derivative $\mathbf{U}_t$ is reshaped into a column vector and the system in Eq.~\eqref{eq1b} is written as:
\begin{equation}
    \mathbf{U}_t = \boldsymbol{\Theta}(\mathbf{U})\boldsymbol{\Xi}.
    \label{PDEFIND}
\end{equation}
For most PDEs, $\boldsymbol{\Xi}$ is sparse and can be identified with a similar sparsity-promoting regression:
\begin{equation}
    \boldsymbol{\Xi} = \argmin_{\Hat{\boldsymbol{\Xi}}} \frac{1}{2} \|\mathbf{U}_t-\boldsymbol{\Theta}(\mathbf{U})\Hat{\boldsymbol{\Xi}}\|_2^2 +  R(\hat{\boldsymbol{\Xi}}). 
\end{equation}
STRidge improves model identification with highly correlated data that is common in PDE regression problems.
PDE-FIND is extremely prone to noise, because noise is amplified when computing high-order partial derivatives for $\boldsymbol{\Theta}$. 
To make PDE-FIND more noise robust, integral~\cite{schaeffer2017sparse} and weak formulations~\cite{reinbold2020using,messenger2021weak2} were introduced.
Instead of discovering a model based on Eq.~\eqref{eq1b}, the PDE can be multiplied by a weight $\mathbf{w}_j(\mathbf{u},t)$ and integrated over a domain $\Omega_k$. This can be repeated for a number of combinations of $\mathbf{w}_j(\mathbf{u},t)$ and $\Omega_k$. Stacking the results of the integration over different domains using different weights leads to a linear system
\begin{equation}
    \mathbf{q}_0 = \mathbf{Q}\boldsymbol{\Xi},
    \label{weak}
\end{equation}
with $\mathbf{q}_0$ and $\mathbf{Q}=[\mathbf{q}_1,\cdots,\mathbf{q}_D]$ the integrated left hand side and integrated library of candidate terms, which replace $\textbf{U}_t$ and the library of nonlinear functions $\boldsymbol{\Theta}(\mathbf{U})$. 
As with PDE-FIND, sparse regression can be employed to identify a sparse matrix of coefficients $\boldsymbol{\Xi}$, using STLS, STRidge, or other regularizers.
For all of our results, we use this weak formulation as a baseline and for the basis of ensemble models.

\section{Ensemble SINDy} \label{method}

In this work, we introduce E-SINDy, which incorporates ensembling techniques into data-driven model discovery. 
Ensembling is a well established machine learning technique
that combines multiple models to improve prediction. 
A range of ensembling methods exist, such as bagging (bootstrap aggregation)~\cite{breiman1996bagging}, bragging (robust bagging)~\cite{buhlmann2003bagging, buhlmann2012bagging}, and boosting~\cite{schapire1990strength, freund1995boosting}.
Structure learning techniques such as cross-validation~\cite{hastie2001elements} or stability selection~\cite{meinshausen2010stability} can also be considered ensembling methods, because they combine and use the information of a collection of learners or models. 
For model discovery, ensembling improves robustness and naturally provides inclusion probabilities and uncertainty estimates for the identified model coefficients, which enable probabilistic forecasting and active learning. 

\begin{figure*}[t]
    \includegraphics[width=\textwidth]{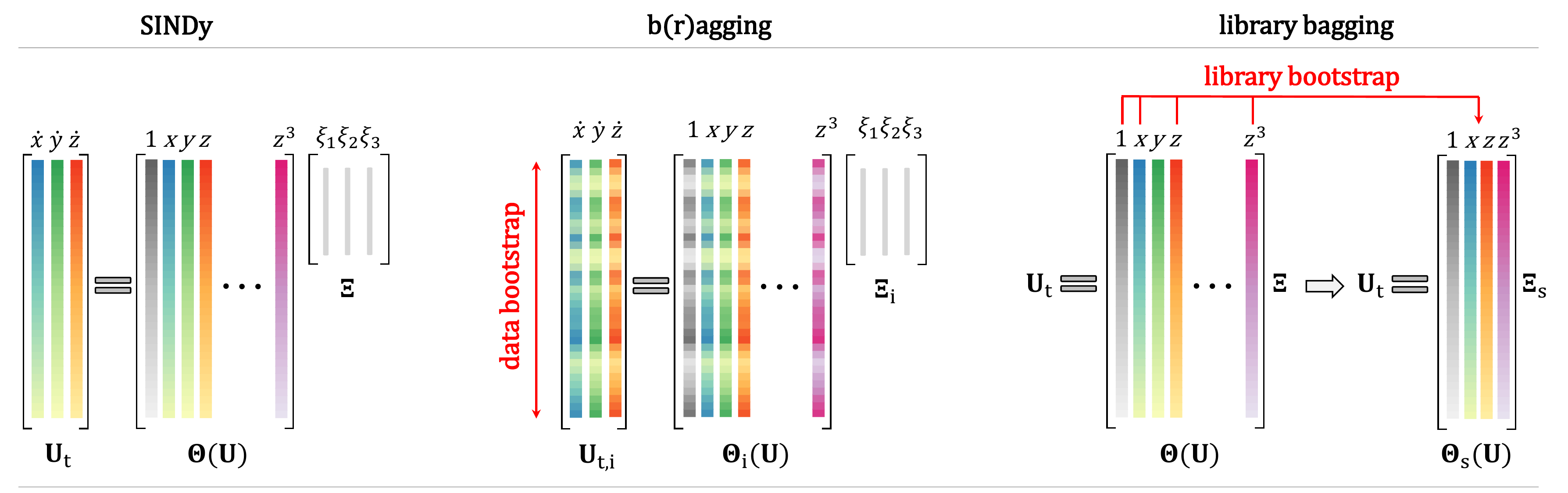}
    \caption{Schematic of SINDy and E-SINDy with b(r)agging and library bagging. Shown is a single model of the ensemble. In case of b(r)agging, data bootstraps (data samples with replacement) are used to generate an ensemble of SINDy models. The E-SINDy model is aggregated by taking the mean of the identified coefficients for bagging, and the median for bragging. In case of library bagging, instead of data bootstraps, library term bootstraps are sampled without replacement. Library terms with low inclusion probability are discarded and the E-SINDy model can be identified on the smaller library using standard SINDy or b(r)agging E-SINDy.}
    \label{fig2}
\end{figure*}

Here, we propose two new ensemble model discovery methods: the first method is called b(r)agging E-SINDy, and the second method is called library E-SINDy. A general schematic of E-SINDy is shown in figure~\ref{fig1}, and a schematic of the sparse regression problems for b(r)agging and library E-SINDy is shown in figure~\ref{fig2}.
Our first method, b(r)agging E-SINDy, uses data bootstraps to discover an ensemble of models that are aggregated by taking the mean of the identified model coefficients in case of bagging, and taking the median in the case of bragging. 
Bootstraps are data samples with replacement. 
Applied to SINDy to identify ODEs, we first build a library of candidate terms $\boldsymbol{\Theta(\textbf{U})}$ and derivatives $\textbf{U}_t$. 
From the $m$ rows of the data matrices $\textbf{U}_t$ and $\boldsymbol{\Theta(\textbf{U})}$, corresponding to $m$ samples in time, we select $q$ bootstrapped data samples and generate $q$ SINDy models in the ensemble.  
For each of these $q$ data bootstraps, $m$ new rows are sampled with replacement from the original $m$ rows of the data matrices.  
On average, each data bootstrap will have around $63\%$ of the entries of the original data matrices, with some of these entries being represented multiple times in the bootstrap; for large $m$ this quantity converges to $1-e^{-1} \approx 0.632$.  
In this way, randomness and subsampling is inherent to the bootstrapping procedure.  
From the $q$ identified SINDy models in the ensemble, we can either directly aggregate the identified models, or first threshold coefficients with low inclusion probability.  
The procedure is illustrated in algorithm~\ref{algoBagging} for bagging E-SINDy using STRidge. The same procedure applies for bragging, taking the median instead of the mean, and using other regularizers than STRidge. 
Note that there are other random data subsampling approaches that may be used, such as generating $q$ models based on $q$ random subsamples of $p<m$ rows of the data without replacement, of which there are ${m \choose p}$.  
However, boostrapping based on selection with replacement is the most standard procedure.  

\begin{lstlisting}[style=interfaces,language=Matlab,caption={Bagging E-SINDy with STLS.},label={algoBagging},mathescape=true]
$\textbf{Inputs}:$  library $\boldsymbol{\Theta}$ and time derivatives $\mathbf{U}_t$
         sparsification and ridge parameters $\lambda_1$ and $\lambda_2$
         number of bootstraps $q$
         inclusion probability threshold $tol$
$\textbf{Outputs}:$ SINDy model coefficients $\boldsymbol{\Xi}$ 
         coefficient inclusion probabilities $ip$
        
function bagging($\boldsymbol{\Theta}$, $\mathbf{U}_t$, $\lambda$, $q$, $tol$)
    % Assemble data matrix that we sample from
    data = [$\mathbf{U}_t$, $\boldsymbol{\Theta}$] 
    % Bootstrap: identifying SINDy models $\boldsymbol{\Xi}_{B}$ (tensor size [m,n,q])
    $\boldsymbol{\Xi}_{B}$ = bootstrap(q,@sparsifyDynamics(data, $\lambda_1$, $\lambda_2$))
    % Compute inclusion probability
    ip = mean($\boldsymbol{\Xi}_{B}$($\boldsymbol{\Xi}_{B}\sim$=0), 3)
    % Aggregate SINDy models
    $\boldsymbol{\Xi}$ = mean($\boldsymbol{\Xi}_{B}$, 3)
    % and threshold coefficients with low inclusion probability
    $\boldsymbol{\Xi}$(ip<tol) = 0
end
\end{lstlisting}
The second method proposed, library bagging E-SINDy, samples library terms instead of data pairs. 
We sample $l$ out of $D$ library terms without replacement. 
In case of sampling library terms, replacement does not affect the sparse regression problem. However, using smaller libraries can drastically speed up model identification, as the complexity of the least squares algorithm is $\mathcal{O}(ml^2)$. Library bagging with small $l$ can therefore help counteract the increased computational cost of solving multiple regression problems in the ensemble. 
As with bagging E-SINDy, we obtain an ensemble of models and model coefficient inclusion probabilities. We can directly aggregate the models and threshold coefficients with low inclusion probabilities to get a library E-SINDy model. We can also use the inclusion probabilities to threshold the library, only keeping relevant terms, and run bagging E-SINDy using the smaller library. This can be particularly useful if we start with a large library: we first identify and remove all library terms that are clearly not relevant and then run bagging E-SINDy on the smaller library.
We show a pseudo code of library bagging E-SINDy in algorithm~\ref{algoLibrary}. 

\begin{lstlisting}[style=interfaces,language=Matlab,caption={Library E-SINDy with STLS.},label={algoLibrary},mathescape=true]
$\textbf{Inputs}:$  library $\boldsymbol{\Theta}$ and time derivatives $\mathbf{U}_t$
         sparsification and ridge parameters $\lambda_1$ and $\lambda_2$
         number of bootstraps $q$
         number of library terms $l$ 
         inclusion probability threshold $tol$
$\textbf{Outputs}:$ SINDy model coefficients $\boldsymbol{\Xi}$ 
         coefficient inclusion probabilities $ip$
        
function libBag($\boldsymbol{\Theta}$, $\mathbf{U}_t$, $\lambda$, $q$, $tol$)
	for b = 1:q 
    		% Sample library
    		$\boldsymbol{\Theta}_b$ = $\boldsymbol{\Theta}$(:,datasample(1:D,l,Replace,false))
    		% Identify SINDy model
    		$\boldsymbol{\Xi}_B$(:,:,b) = @sparsifyDynamics($\boldsymbol{\Theta}_b$, $\mathbf{U}_t$, $\lambda_1$, $\lambda_2$)
	end
	% Compute inclusion probability
	ip = mean($\boldsymbol{\Xi}_B$($\boldsymbol{\Xi}_B\sim$=0), 3)
	% Aggregate SINDy model
	$\boldsymbol{\Xi}$ = mean($\boldsymbol{\Xi}_{B}$, 3)
	% Threshold coefficients with low inclusion probability
	$\boldsymbol{\Xi}$(ip<tol) = 0
	% or use thresholded library $\boldsymbol{\Theta}_{t}$ for bagging
	$\boldsymbol{\Theta}_{t}$ = $\boldsymbol{\Theta}$($\boldsymbol{\Xi}$(ip<tol))
	$\boldsymbol{\Xi}$ = @bagging($\boldsymbol{\Theta}_{t}$, $\mathbf{U}_t$, $\lambda$, $q$, $tol$)
end
\end{lstlisting}
E-SINDy provides inclusion probabilities and uncertainty estimates for the discovered model coefficients, thus connecting to Bayesian model identification techniques. 
The identified ensemble of model coefficients can be used to compute coefficient probability density functions, which form a posterior distribution $p(\boldsymbol{\Xi}|\boldsymbol{X})$. 
In terms of forecasting, we can either use the aggregated mean or median of the identified coefficients to forecast, or we can draw from multiple identified SINDy models to generate ensemble forecasts that represent posterior predictive distributions $p(\boldsymbol{x}(t)|\boldsymbol{X})$ that provide prediction confidence intervals.

\section{Results} \label{results}

We now apply E-SINDy to challenging synthetic and real world data sets to identify ordinary and partial differential equations. We apply library bagging E-SINDy to a real world ecological data set, showing its performance in the very low data limit. For PDEs, we use the recent weak-SINDy~\cite{messenger2021weak2} as a baseline and show the improved noise robustness when using E-SINDy for identifying a range of PDEs. 
Trends for the noise and data length sensitivity of bagging, bragging, and library bagging to identify the chaotic Lorenz system dynamics are presented in the appendix. 

\subsection{Ordinary differential equations} \label{ODE}

\begin{figure*}[t]
    \includegraphics[width=\textwidth]{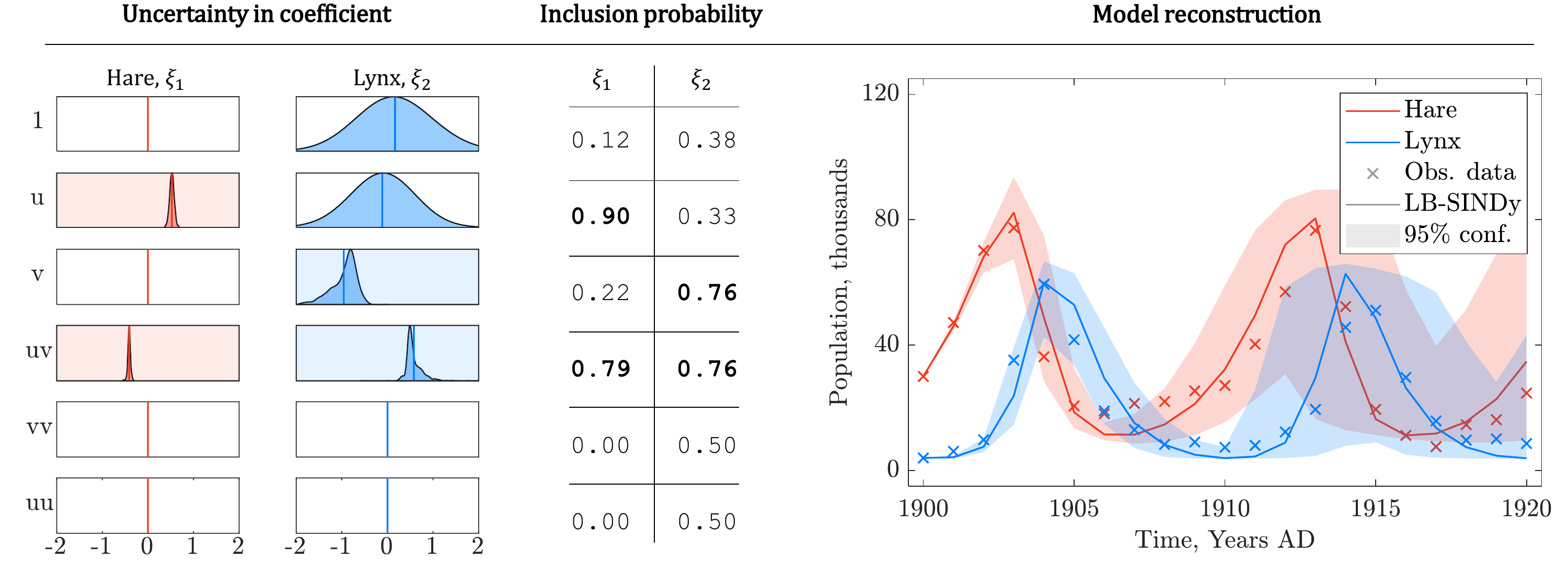}
    \caption{Library bagging E-SINDy (LB-SINDy) on real data: data consisting of measurements by the Hudson Bay Company of lynx and hare pelts between 1900 and 1920. (Left) Uncertainty in the identified model coefficients, (center) inclusion probabilities of the model coefficients, and (right) model reconstruction. LB-SINDy (continuous lines) uses the mean value of the identified coefficients for reconstruction, and the $95\%$ confidence interval depicts ensemble reconstruction, drawing five models and averaging the coefficients for 1000 realizations.}
    \label{fig4}
\end{figure*}

We apply E-SINDy to a challenging real-world data set from the Hudson Bay Company, which consists of the yearly number of lynx and hare pelts collected between 1900 and 1920.  
These pelt counts are thought to be roughly proportional to the population of the two species~\cite{hewitt1921conservation}. 
Lynx are predators whose diet depends on hares. The population dynamics of the two species should, therefore, be well approximated by a Lotka-Volterra model.
There are several challenges in identifying a SINDy model from this data set:  there are only 21 data points available, and there is large uncertainty in the measurements arising from weather variability, consistency in trapping, and other changing factors over the years measured.
In figure~\ref{fig4}, we show that E-SINDy correctly identifies the Lotka-Volterra dynamics, providing model coefficient and inclusion probabilities and confidence intervals for the reconstructed dynamics. 
We use library bagging followed by bagging to identify a sparse model in this very low data limit with only 21 data points per species. 
Similar results for the lynx-hare data set were recently published using a probabilistic model discovery method~\cite{hirsh2021sparsifying} based on sparse Bayesian inference.  
This approach employed Markov Chain Monte Carlo, for which the computational effort to generate a probabilistic model is comparably high, taking several hours of CPU time. 
In contrast, E-SINDy takes only seconds to identify a model and its coefficient and inclusion probabilities. 

\subsection{Partial differential equations}
\label{PDE}

In this section, we present results applying E-SINDy to discover partial differential equations from noisy data.
We use the recent weak-SINDy (WSINDy) implementation~\cite{messenger2021weak2} as the baseline model for ensembling. WSINDy was successfully applied to identify models in the high-noise regime using large libraries. 
We perform library bagging on the system of Eq.~\ref{weak} instead of Eq.~\ref{PDEFIND}, and refer to the resulting method as ensemble weak SINDy (E-WSINDy). 

\begin{figure*}[t]
    \includegraphics[width=\textwidth]{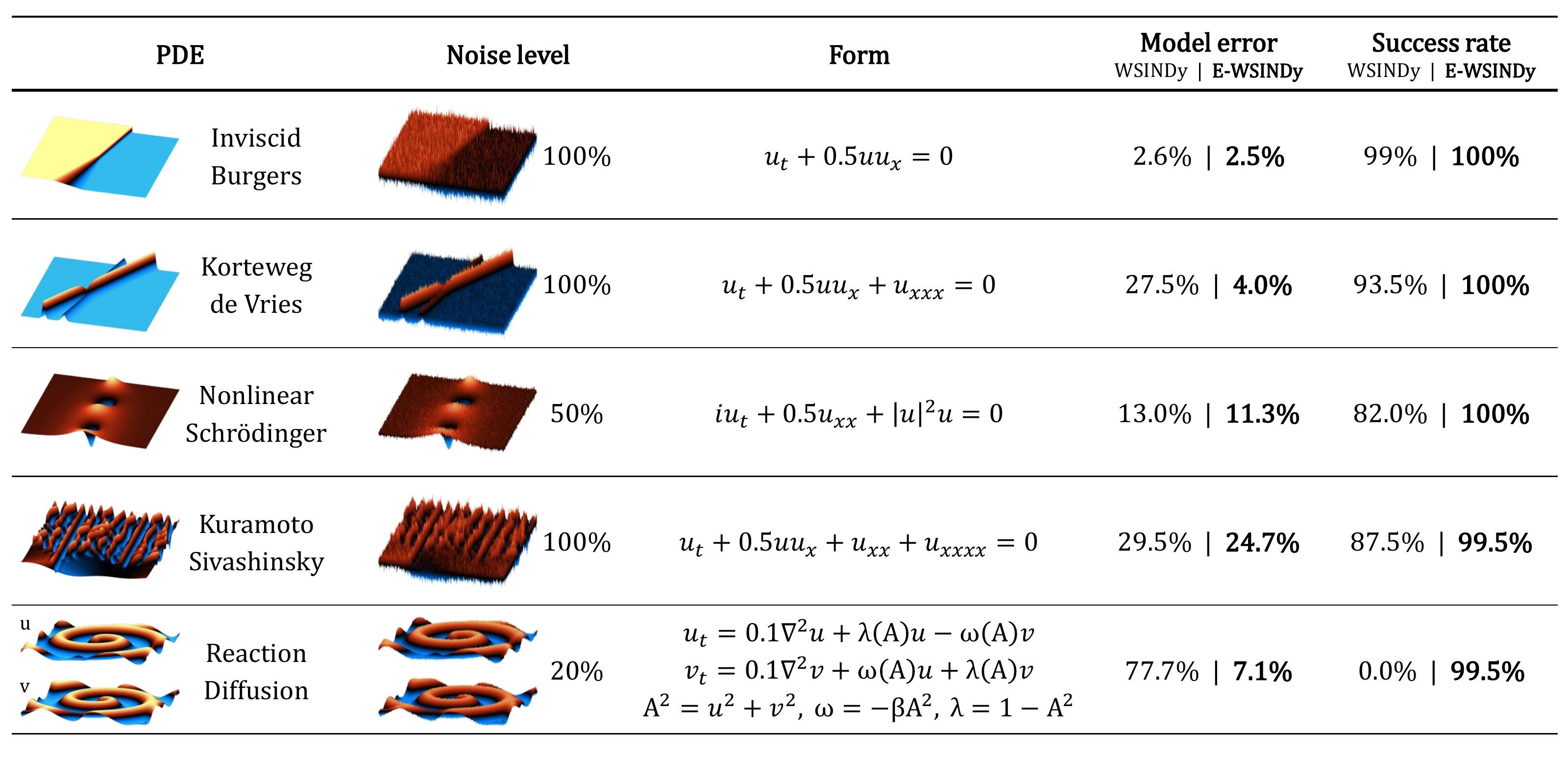}
    \vspace{-.35in}
    \caption{Comparison of model error and success rate of discovered PDEs using weak-SINDy and ensemble weak-SINDy. Ensembling robustifies and improves the accuracy of PDE identification.}
    \label{fig5}
\end{figure*}
We apply E-WSINDy to identify PDEs from synthetic data for the inviscid Burgers, Korteweg de Vries, nonlinear Schroedinger, Kuramoto Sivashinsky, and reaction diffusion equations. 
We quantify the accuracy and robustness of the model identification by assessing the success rate and model coefficient errors for a number of noise realizations. 
The success rate is defined as the rate of identifying the correct non-zero and zero terms in the library, averaged over all realizations.
The model coefficient error $E_{c}$ quantifies how much the identified coefficients $\hat{\boldsymbol{\Xi}}$ deviate from the true parameters $\boldsymbol{\Xi}$ that we use to generate the data: 
\begin{equation}
    E_{c}=\frac{\|\boldsymbol{\Xi}-\hat{\boldsymbol{\Xi}}\|_{2}}{\|\boldsymbol{\Xi}\|_{2}}.
\end{equation}
The results are summarized in figure~\ref{fig5}. 
For all PDEs, E-WSINDy reduces the model coefficient error and increases the success rate of the model discovery. 
Moreover, E-WSINDy can accurately identify the correct model structure for the reaction diffusion case, where WSINDy falsely identifies a linear oscillator model instead of the nonlinear reaction diffusion model.
To investigate the limits of E-WSINDy, we further increase the noise level for each case up to the point where the success rate drops below $90\%$. 
On average, for all investigated PDEs, we find that ensembling improves noise robustness of WSINDy by a factor of 2.3.
We conclude that ensembling significantly improves model discovery robustness and enables the identification of PDEs in the extreme noise limit.

\subsection{Exploiting ensemble statistics for active learning}
\label{AL}

\begin{figure*}[t]
    \includegraphics[width=\textwidth]{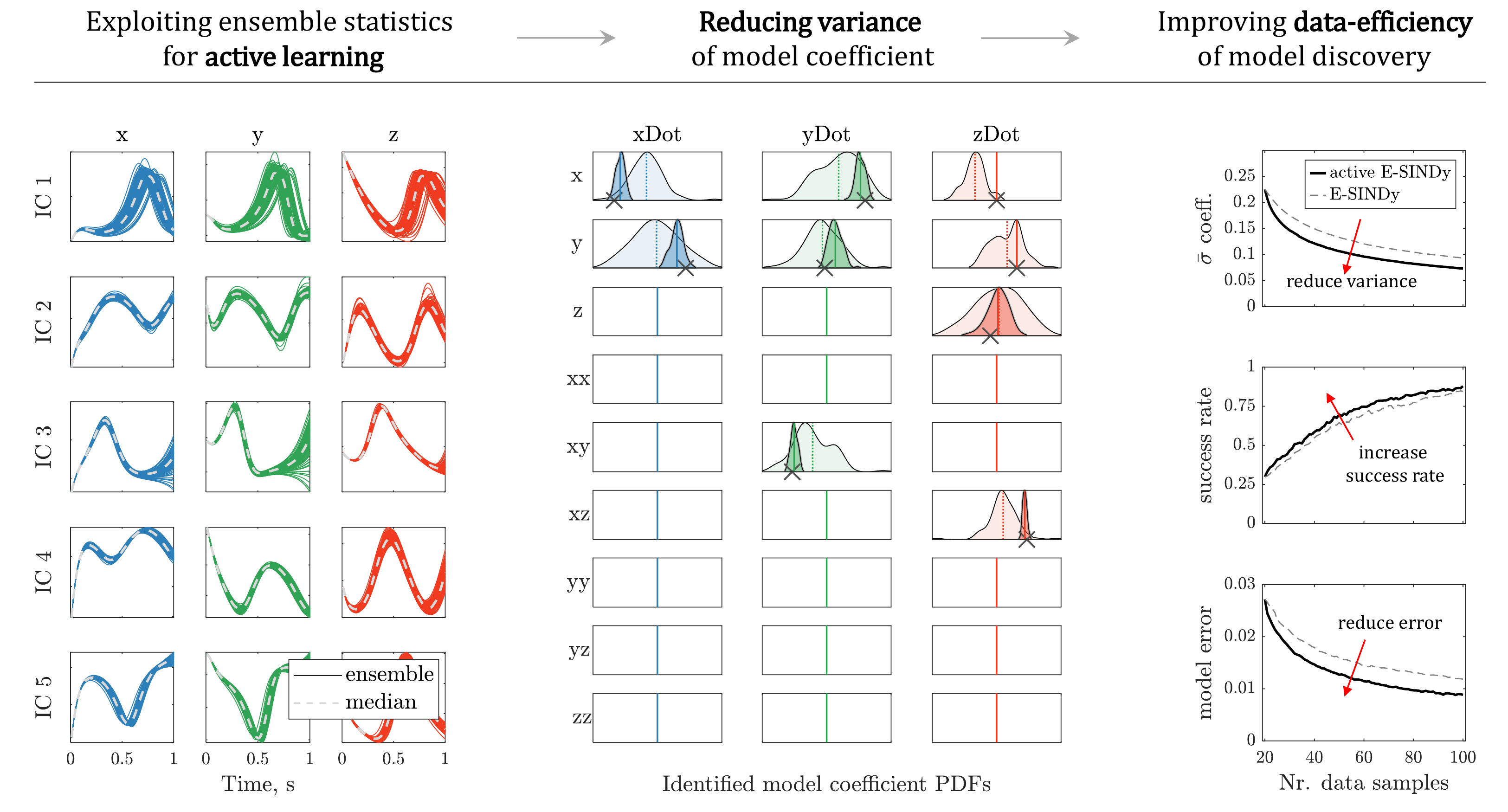}
    \caption{Exploiting ensemble statistics for active learning. Active E-SINDy randomly selects a number of initial conditions (IC), computes the ensemble forecast variance at each IC, and explores the IC with highest variance. Left: ensemble forecasts from different ICs. Center: reduced variance of identified model coefficients after several active learning steps. Right: Improved data-efficiency and accuracy of model discovery using active learning.}
    \label{fig6}
\end{figure*}

We now present results exploiting the ensemble statistics for active learning~\cite{krogh1995neural,zhu2010active}. 
Active learning is a machine learning method that can reduce training data size while maintaining accuracy by actively exploring regions of the feature space that maximally inform the learning process~\cite{settles2009active, settles2011theories}.
Here, we leverage the ensemble statistics of E-SINDy to identify and sample high-uncertainty regions of phase space that maximally inform the sparse regression.
In E-SINDy, we collect data from a single initial condition or from multiple randomly selected initial conditions and identify a model in one shot. 
Instead, we can successively identify E-SINDy models and exploit their ensemble statistics to identify new initial conditions with high information content, which can improve the data-efficiency of the model discovery process. 
The basic idea is to compute ensemble forecasts from a large set of initial conditions using successively improved E-SINDy models and only explore regions with high ensemble forecast variance. 
Our simple but effective active E-SINDy approach iteratively identifies models in three steps:
(1) collecting a small amount of randomly selected data to identify an initial E-SINDy model;  
(2) selecting a number of random initial conditions and computing the ensemble forecast variance for each initial condition using the current E-SINDy model; and 
(3) sampling the true system with the initial condition with highest variance.
Finally, we concatenate the newly explored data to the existing data set to identify a new E-SINDy model, and continue the model identification until the model accuracy and/or variance of the identified model coefficients converge.

Here, we test active E-SINDy on the Lorenz system dynamics introduced in figure~\ref{fig1} and appendix~\ref{app1} and show the results in figure~\ref{fig6}.
On the left, we show five illustrative ensemble forecasts from different initial conditions after initializing the algorithm. In total, at each iteration, we compute ensemble forecasts at 200 different initial conditions. We found that at each initial condition, integrating a single time step forward in time is informative enough to compute ensemble forecasting variance. In the center of figure~\ref{fig6}, we show the probability density functions of the identified model coefficients at initialization have wide distributions, and after 80 active learning steps the variance of the distributions is significantly reduced. 
 Figure~\ref{fig6} (right) also shows the improved data-efficiency of the model discovery using active learning E-SINDy compared to E-SINDy. Through active E-SINDy, we reduce the variance of the identified model coefficients, increase the success rate of identifying the correct model structure, and reduce the model coefficient error compared to standard E-SINDy.

\subsection{Ensemble SINDy model predictive control}

\begin{figure*}[t]
    \includegraphics[width=\textwidth]{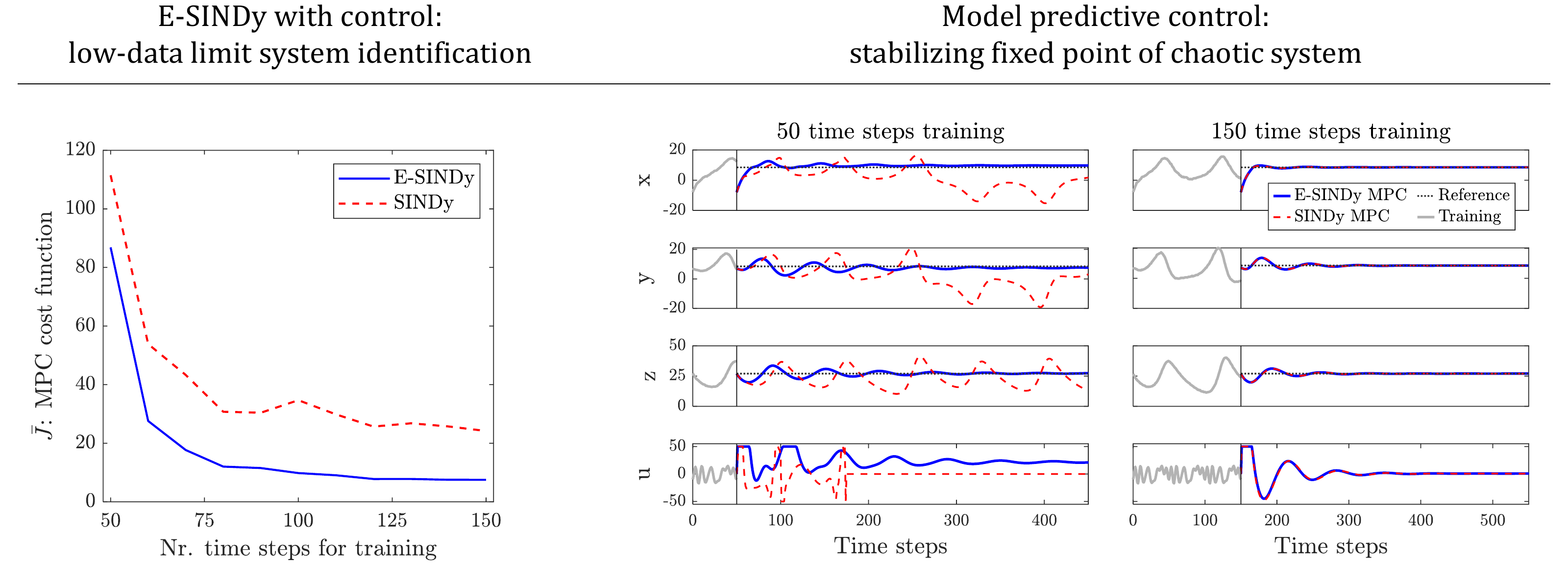}
    \caption{System identification in the low-data limit for model predictive control. Left: MPC cost function average $\bar{J}$ over number of time steps used for training with E-SINDy (blue continuous line) and SINDy (red dashed line). Right: controlled trajectory coordinates x, y, z of Lorenz system with model predictive control input u, for models trained with E-SINDy (blue continuous line) and SINDy (red dashed line) using 50 and 150 time step data points.}
    \label{fig7}
\end{figure*}

It is also possible to use E-SINDy to improve model predictive control (MPC)~\cite{garcia1989model,morari1999model,mayne2014automatica}.
MPC is a particularly compelling approach that enables control of strongly nonlinear systems with constraints, multiple operating conditions, and time delays. 
The major challenge of MPC lies in the development of a suitable model. %
Deep neural network models have been increasingly used for deep MPC~\cite{Peng.2009,Zhang2016icra}; however, they often rely on access to massive data sets, have limited ability to generalize, do not readily incorporate known physical constraints, and are computationally expensive.
Kaiser et al.~\cite{Kaiser2018prsa} showed that sparse models obtained via SINDy perform nearly as well with MPC, and may be trained with relatively limited data compared to a neural network. 
Here we show that E-SINDy can further reduce the training data requirements compared to SINDy, enabling the control of nonlinear systems in the very low data limit.

We use E-SINDy to identify a model of the forced Lorenz system dynamics and use MPC to stabilize one of the two unstable fixed points $(\pm\sqrt{72},\pm\sqrt{72},27)$.
The Lorenz system is introduced in figure~\ref{fig1} and we add a control input $u$ to the first state of the dynamics: $\dot{x}=\sigma(y-x)+u$. 
The control problem is based on Kaiser et al.~\cite{Kaiser2018prsa}, where the MPC cost function and weight matrices are described in detail.
In figure~\ref{fig7}, we show the performance of MPC based on E-SINDy models for different training data length and $\text{noise}=0.01$. On the left, we show the sensitivity of the mean MPC cost to training data length. We run 1000 noise realizations and average the mean MPC cost of all runs. The right panel shows trajectories of the controlled Lorenz system for models trained with E-SINDy and SINDy, using 50 and 150 time step data points.
E-SINDy significantly improves the MPC performance in the low data limit compared to SINDy.

\section{Discussion}

This work has developed and demonstrated a robust variant of the SINDy algorithm based on ensembling.  
The proposed E-SINDy algorithm significantly improves the robustness and accuracy of SINDy for model discovery, reducing the data requirements and increasing noise tolerance.  
E-SINDy exploits foundational statistical methods, such as bootstrap aggregating, to identify ensembles of ordinary and partial differential equations that govern the dynamics from noisy data. 
From this ensemble of models, aggregate model statistics are used to generate inclusion probabilities of candidate functions, which promotes interpretability in model selection and provides probabilistic forecasts. 
We show that ensembling may be used to improve several standard SINDy variants, including the integral formulation for PDEs.  
Combining ensembling with the integral formulation of SINDy enables the identification of PDE models from data with more than twice as much measurement noise as has been previously reported.  
These results are promising for the discovery of governing equations for complex systems in neuroscience, power grids, epidemiology, finance, or ecology, where governing equations have remained elusive.
Importantly, the computational effort to generate probabilistic models using E-SINDy is low.
E-SINDy produces accurate probabilistic models in seconds, compared to existing Bayesian inference methods that take several hours.
Library bagging has the additional advantage of making the least squares computation more efficient by sampling only small subsets of the library.
E-SINDy has also been incorporated into the open-source PySINDy package~\cite{deSilva2020JOSS,kaptanoglu2021pysindy} to promote reproducible research.

We also present results exploiting the ensemble statistics for active learning and control.
Recent active exploration methods~\cite{shyam2019model} and active learning of nonlinear system identification~\cite{mania2020active} suggest exploration techniques using trajectory planning to efficiently explore high uncertainty regions of the feature space. 
We use the uncertainty estimates of E-SINDy to explore high uncertainty regions that maximally inform the learning process. 
Active E-SINDy reduces the variance of the identified model coefficients, increases the success rate of identifying the correct model structure, and reduces the model coefficient error compared to standard E-SINDy in the low data limit.
Finally, we apply E-SINDy to improve nonlinear model predictive control.
SINDy was recently used to generate models for real-time MPC of nonlinear systems. We show that E-SINDy can significantly reduce the training data required to identify a model, thus enabling control of nonlinear systems with constraints in the very low data limit.
An exciting future extension of the computationally efficient probabilistic model discovery is to combine the active learning and MPC strategies based on E-SINDy. 
Highly efficient exploration and identification of nonlinear models may also enable learning task-agnostic models that are fundamental components of model-based reinforcement learning.

\section*{Acknowledgments}

The authors acknowledge support from the Air Force Office of Scientific Research (AFOSR FA9550-19-1-0386), the Army Research Office (ARO W911NF-19-1-0045), and the National Science Foundation AI Institute in Dynamic Systems (Grant No. 2112085).

\appendix
\section{E-SINDy trends compared with standard SINDy}
\label{app1}

Here, we evaluate noise and data length sensitivity of E-SINDy compared to the standard SINDy implementation with no innovations for noise robustness.  
This case does not demonstrate the maximum noise robustness, but rather provides a simple and intuitive introduction to how E-SINDy compares with a baseline SINDy algorithm.  
E-SINDy can be used with most SINDy variants, and similar performance increases are found.  
For example, using integral SINDy or computing derivatives using the total variation derivative will dramatically improve the baseline robustness and therefore the E-SINDy robustness, as shown in section~\ref{PDE}.  

\begin{figure*}[t]
    \includegraphics[width=\textwidth]{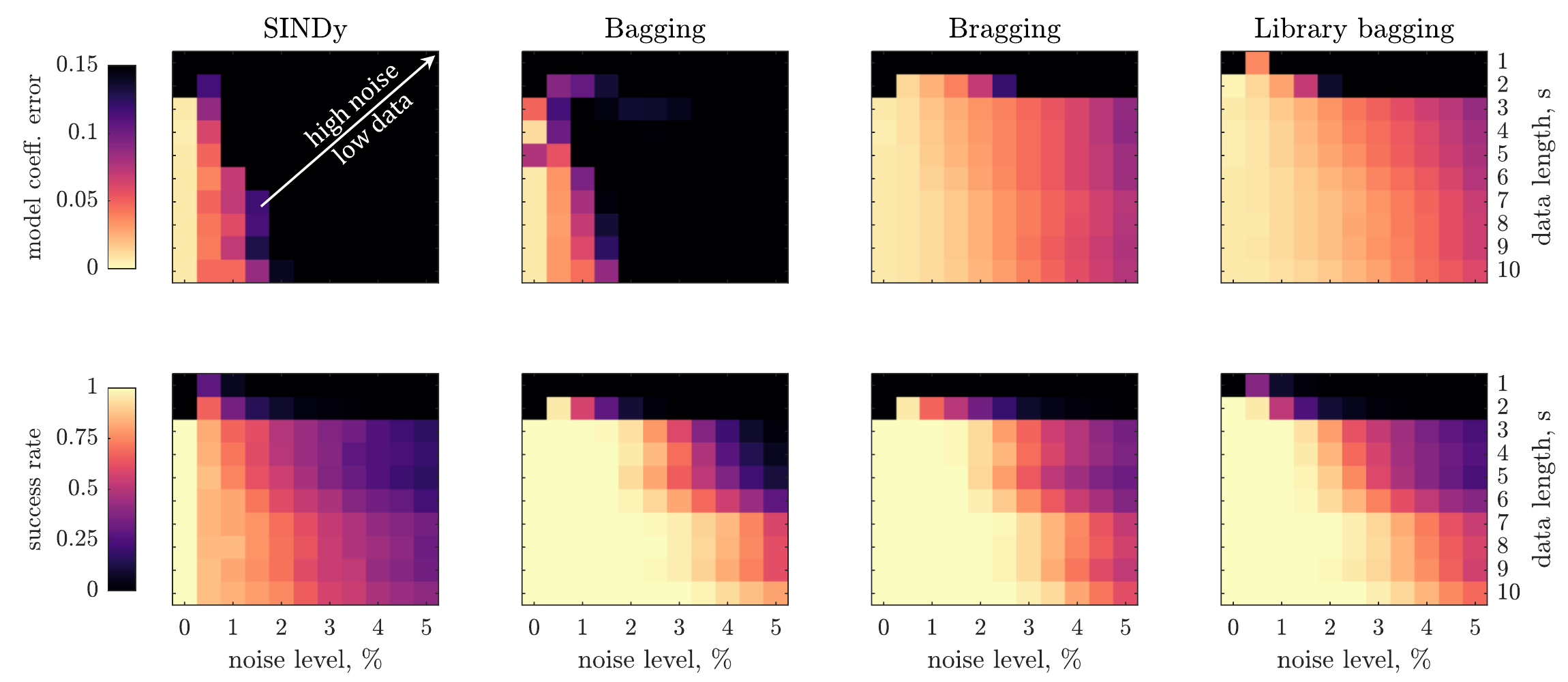}
    \vspace{-.3in}
    \caption{Noise level and data length sensitivity of different SINDy and ensemble SINDy  models of the Lorenz system dynamics. (Top) Model coefficient error, and (bottom) success rate in identifying the correct model structure. Each noise level and data length case is run on 1000 realizations of noise and averaged for plotting. Ensembling improves the accuracy of SINDy with less data and more noise.}
    \label{fig3}
\end{figure*}

We define three commonly used metrics of success for identifying sparse nonlinear models: (1) the structure of the identified model, i.e. identifying the correct non-zero terms in library; (2) the error of the identified model coefficients; and (3) the accuracy of predicted trajectories. 
The accuracy of a predicted trajectory is highly sensitive to small differences in the model coefficients, especially for chaotic systems. Therefore, we assess the structure of the identified model and the error of the identified coefficients quantitatively, and assess the accuracy of the forecast qualitatively. 
In figure~\ref{fig3}, we compare the noise level and data length sensitivity of SINDy with bagging, bragging, and library bagging E-SINDy, for data from the Lorenz system. 
We generate trajectories from the Lorenz system with system parameters shown in figure~\ref{fig1} and initial condition $\textbf{u}_0=[-8,\hspace{1pt} 7,\hspace{1pt} 27]$, for a range of data length. 
We then add Gaussian white noise with zero mean and variance $ {\sigma}/{\left\|\mathbf{u}\right\|_{\text{rms}}}$
where $\left\|\mathbf{u}\right\|_{\text{rms}}$ is the root mean squared value of the trajectory data.
We test each noise level and data length case with 1000 random realizations of noise, compute the success rate of identifying the correct model structure in $\boldsymbol{\Xi}$, and calculate the mean error of the identified model coefficients. 
The sparsity promoting hyperparameter is set to $\lambda = 0.2$ for all methods. The number of models in the ensemble is set to $q=100$. The inclusion probability used to threshold model coefficients is set to $tol=0.6$ for b(r)agging and $tol=0.4$ for library bagging. We choose a smaller threshold for library bagging because we perform bagging after library bagging. Thus, the goal is not to aggressively threshold the identified coefficients, but only to reduce the initial library size to simplify the complexity of the subsequent bagging E-SINDy. 
We found that this strategy can improve the success rate and reduce the model error. Additionally, in the case of library bagging, we define the number of library terms in each sample to be $l=0.6D$.

The results for the different ensemble methods over noise level and data length are shown in figure~\ref{fig3}. We see that ensembling clearly improves the accuracy and robustness of the model discovery, enabling model identification with less data and more noise than standard SINDy. 
As expected, bragging further robustifies bagging, and library bagging also improves accuracy and robustness. 
Additionally, library bagging has the advantage that we can consider larger libraries, discarding irrelevant terms in the library before applying b(r)agging E-SINDy.

In terms of the accuracy of the E-SINDy forecast, we show library bagging for $2.5\%$ noise and data length $T=10s$ in the bottom of figure~\ref{fig1}.
Library bagging E-SINDy (colored solid line) is able to forecast several chaotic lobe switches, compared to SINDy (gray dotted line), which deviates much earlier from the true dynamics (dashed line). 
We also show the $95\%$ confidence intervals for different ensemble forecasting strategies, where multiple models are drawn from the bootstrap and used for prediction.  
Additionally, figure~\ref{fig1} shows the model coefficient distribution and inclusion probabilities of the identified library bagging E-SINDy model. 
The variance of the identified coefficients is comparably low, and the inclusion probabilities for the non-zero terms are clearly higher than for the zero terms. 
    
As an alternative to bagging and bragging, we also perform the same sensitivity analysis to noise level and data length using the stability selection method~\cite{meinshausen2010stability} that was recently used for PDE model discovery~\cite{maddu2019stability}. Stability selection can be used to automatically determine the level of regularization in SINDy. The method randomly samples data subsets of size $m/2$ without replacement, identifies SINDy models for each subset, computes an importance measure and thresholds library terms, and solves a final linear least-squares regression on the thresholded library. The method in its original form performs slightly worse than b(r)agging E-SINDy. However, the success rate of stability selection SINDy can be improved, achieving comparable performance to bragging E-SINDy, by sampling subsets of size $m$ with replacement instead of $m/2$ without replacement.

\bibliographystyle{unsrt}
 \begin{spacing}{.9}
 \small{
 \setlength{\bibsep}{2.9pt}
 \bibliography{references}
 }
 \end{spacing}
\end{document}